%% file: 0.tex
\documentclass[11pt]{amsart}
\usepackage{graphicx, epstopdf, verbatim, fancyhdr, indent, wrapfig}
\usepackage{subfigure}

\pdfoutput=1

\usepackage[T1]{fontenc}
\newcommand{\scb}[1]{\textbf{\textsc{#1}}}
 
\usepackage[nohug]{diagrams}\diagramstyle[labelstyle=\scriptstyle]
\usepackage[usenames,dvipsnames]{color}
\usepackage[colorlinks=true, urlcolor=MidnightBlue, linkcolor=MidnightBlue, citecolor=MidnightBlue, pdftitle={Generalized Gauss maps and integrals for three-component links: Toward higher helicities for magnetic fields and fluid flows, Part II}, pdfauthor={Dennis DeTurck, Herman Gluck, Rafal Komendarczyk, Paul Melvin, Haggai Nuchi, Clayton Shonkwiler and David Shea Vela-Vick}, 
pdfsubject={Geometric Topology, Knot Theory}, pdfkeywords={mu invariant, Hopf invariant, link-homotopy, Gauss map, helicity}]{hyperref}

\usepackage[left=1.3in,top=1.3in,right=1.3in,bottom=1.2in,head=.1in]{geometry}

\bfseries
\def\part#1{\smallbreak\noindent{\textbf{\boldmath #1}}\smallbreak} %\textsc
\newtheorem{theorem}{Theorem}

\newtheorem*{thmBcont}{Theorem B (continued)}
\newtheorem*{bridge}{Bridge Lemma}
\newtheorem{proposition}{Proposition}

\theoremstyle{definition}

\newcommand{\thmref}[1]{Theorem~\ref{thm:#1}}

\def\|{|}

\newcommand{\itb}{\bfseries\itshape\boldmath}

\def\bz{\mathbb Z}

\def\br{\mathbb R}

\newarrow{Equals}=====

\def\x{\textup{\bf x}}

\def\a{\textup{\bf a}}

\def\b{\textup{\bf b}}
\def\c{\textup{\bf c}}

\def\n{\textup{\bf n}}

\def\bah#1{\overline#1}

\newcommand{\reals}{\br}
\newcommand{\ints}{\bz}
\newcommand{\zero}{{\bf0}}

\def\Vert{|}
\newcommand{\dast}{*} %{{\displaystyle{\ast}}}
\newcommand{\Hopf}{\mathop{\rm Hopf}\nolimits}

\newcommand{\Isom}{\mathop{\rm Isom}\nolimits}
\renewcommand{\Re}{\mathop{\rm Re}}
\renewcommand{\Im}{\mathop{\rm Im}}

\def\dt{\,\raisebox{-.6ex}{\huge$\hskip-.02in\cdot$\hskip.005in}}
\def\dte{\,\raisebox{-.3ex}{\Large$\hskip-.02in\cdot$\hskip.005in}}

\def\dtw{\,\dt\,}

\def\con#1#2{\mathrm{Conf}_{#1}#2}
\def\conf{\con{3}{\br^3}}
\def\confs{\con{3}{S^3}}

\def\lk{\mathrm{Lk}}

\definecolor{myred}{RGB}{230,36,15}

\begin{document}
%%%%%%%%%%%%%%%

%%%%%%%%%%%%%%%%
%% TITLE and ABSTRACT %%
%%%%%%%%%%%%%%%%

%\vskip-.2in 
% \vskip-.2in

\newgeometry{left=1.3in,top=0.5in,right=1.3in,bottom=1.2in,head=.1in}

\title{Generalized Gauss maps and integrals \\ for three-component links: \\
{\footnotesize toward higher helicities for magnetic fields and fluid flows \\
Part 2}}

%\author{Dennis DeTurck \and Herman Gluck \and Rafal Komendarczyk \and \\ Paul Melvin \and Clayton Shonkwiler \and David Shea Vela-Vick}

%\begin{abstract}
\maketitle

\vskip-.2in

\centerline{Dennis DeTurck, Herman Gluck, Rafal Komendarczyk}\centerline{Paul Melvin, Haggai Nuchi, Clayton Shonkwiler and David Shea Vela-Vick}

\vskip.2in

\centerline{\scb{Prologue}}

\vskip-.1in

{\small
\parskip 3pt

\ 

{\scb{Background.}} The helicity of a magnetic field or fluid flow measures the extent to which its orbits wrap and coil around one another. Woltjer introduced this notion in the late 1950s during his study of the magnetic field in the Crab Nebula, showed that the helicity remains constant as the field evolves according to the equations of ideal magnetohydrodynamics, derived from this a lower bound for the changing field energy, and calculated the stable field at the end of the evolution. The term ``helicity'' was coined ten years later by Moffatt, who rewrote Woltjer's integral formula to reveal its analogy with Gauss's linking integral for two disjoint closed curves in 3-space.

When the helicity of a magnetic field or fluid flow is zero, the lower bound for energy that it provides is also zero, and one hopes for a ``higher order helicity'' which can provide its own lower bound for energy. Monastyrsky and Retakh \cite{MonastyrskyRetakh} and Berger \cite{Berger90} prepared the way for this via integral formulas derived from the Massey product formulation of Milnor's triple linking number of a three-component link, and various authors since then have shown in special cases how to derive nonzero lower energy bounds.

{\scb{What we do here.}} We describe a new approach to triple linking invariants and integrals, aiming for a simpler, wider and more natural applicability to the search for higher order helicities.

To each three-component link in Euclidean 3-space, we associate a {\itb generalized Gauss map} from the 3-torus to the 2-sphere, and show that the pairwise linking numbers and Milnor triple linking number that classify the link up to link homotopy correspond to the Pontryagin invariants that classify its generalized Gauss map up to homotopy.

When the pairwise linking numbers are all zero, we give an integral formula for the triple linking number analogous to the Gauss integral for the pairwise linking numbers, but patterned after J.H.C.\ Whitehead's integral formula for the Hopf invariant, and hence interpretable as the ordinary helicity of a related vector field on the 3-torus.

{\scb{What's new about this?}} Our generalized Gauss map from the 3-torus to the 2-sphere is a natural extension of Gauss's original map from the 2-torus to the 2-sphere; like its predecessor it is equivariant with respect to orientation-preserving isometries of the ambient space, attesting to its naturality and positioning it for application to physical situations; it applies to all three-component links, not just those with pairwise linking numbers zero; and when the pairwise linking numbers are zero, it provides a simple and direct integral formula for the triple linking number which is a natural successor to the classical Gauss integral, with an integrand invariant under orientation-preserving isometries of the ambient space.

{\scb{Application.}} Komendarczyk \cite{Komendarczyk,Komendarczyk10} has applied this approach in special cases to derive a higher order helicity for magnetic fields whose ordinary helicity is zero, and to obtain from this nonzero lower bounds for the field energy.

In the first paper of this series \cite{DGKMSV2}, hereafter ``Part 1'', we did all of the above for three-component links in the three-sphere. The first step there was to find a geometrically natural generalized Gauss map. That same first step is taken here in Euclidean 3-space, but the map itself is entirely different because the requirement of geometric naturality involves a different, and in this case non-compact, group of isometries. After describing this new version of the generalized Gauss map, we build a bridge between the spherical and Euclidean versions, across which we transport proofs and save labor.

\restoregeometry

\parskip 4pt
\addtolength{\baselineskip}{.2pt}

\input{1introduction}

\input{2thmA}

\input{3thmB}

\input{4epilogue}

\input{5bibliography}

\end{document}

%% file: 1introduction.tex
%%%%%%%%%%%%%%%%%
%% SECTION 1: Introduction %%
%%%%%%%%%%%%%%%%%

\section{Introduction}
\label{sec:intro}

%%%%%%%%%%%%%%%%%%%%%%%%%%%%%
%%%%%%%%%%%%%%%%%%%%%%%%%%%%%
\part{Setting the stage.}
%%%%%%%%%%%%%%%%%%%%%%%%%%%%%
%%%%%%%%%%%%%%%%%%%%%%%%%%%%%

Three-component links in Euclidean 3-space $\reals^3$	were classified up to 
{\itb link homotopy}  -- a deformation during which each component may cross itself but distinct components must remain disjoint -- by John Milnor in his senior thesis, published in 1954. A complete set of invariants is given by the pairwise linking numbers $p$ , $q$ and $r$ of the components, and by the \textit{\textbf{triple linking number}}\rm, which is the residue class $\mu$ of one further integer modulo the greatest common divisor of $p$, $q$ and $r$.

For example, the Borromean rings shown below have $p = q = r = 0$ and $\mu=\pm 1$, where the sign depends on the ordering and orientation of the components.

%%%% FIGURE 0: Walnut Borromean rings %%%%%
\begin{figure}[h!]
	\begin{center}
		\includegraphics[height=120pt]{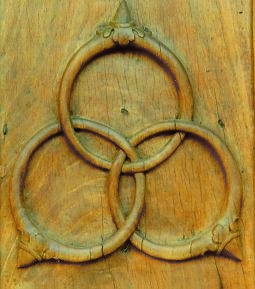}\\
		{\small \bf Borromean Rings} \\
		{\small This is a photograph, courtesy of Peter Cromwell, of a panel in the carved walnut doors of the Church of San Sigismondo in Cremona, Italy.}
	\end{center}
\end{figure}
%%%%%%%%%%%%%%%%%%%

To each ordered, oriented three-component link $L$ in $\reals^3$, we will associate a {\itb generalized Gauss map} $g_L$ from the 3-torus $T^3= S^1 \times S^1 \times S^1$ to the 2-sphere $S^2$, in such a way that link homotopies of $L$ become homotopies of $g_L$. The definition of $g_L$ will be given below.

Maps from $T^3$ to $S^2$ were classified up to homotopy by Lev Pontryagin in 1941. A complete set of invariants is given by the degrees $p$, $q$ and $r$ of the restrictions to the 2-dimensional coordinate subtori, and by the residue class $\nu$ of one further integer modulo twice the greatest common divisor of $p$ , $q$ and $r$, the {\itb Pontryagin invariant} of the map.

This invariant is an analogue of the Hopf invariant for maps from $S^3$ to $S^2$, and is an absolute version of the relative invariant originally defined by Pontryagin for \emph{pairs} of maps from a 3-complex to the 2-sphere that agree on the 2-skeleton of the domain.

\enlargethispage*{1.2em}

Our first main result, \thmref{A} below, equates Milnor's and Pontryagin's invariants $p$, $q$ and $r$ for $L$ and $g_L$, and asserts that
$$2\mu(L) = \nu(g_L).$$
In the special case when $p = q = r = 0$, we derive an explicit and geometrically natural integral formula for the triple linking number, generalizing Gauss's classical integral formula for the pairwise linking number and patterned after J.H.C.\ Whitehead's integral formula for the Hopf invariant. This formula and variations of it are presented in \thmref{B} below.

In the rest of this introduction, we give the background and motivation for our work, then lead up to and provide the definition of the generalized Gauss map of a three-component link in $\reals^3$, give careful statements of Theorems~\ref{thm:A} and~\ref{thm:B}, and finally present the results of a numerical calculation of Milnor's triple linking number using \thmref{B}.

%%%%%%%%%%%%%%%%%%%%%%%%%%%%%%%%%%
%%%%%%%%%%%%%%%%%%%%%%%%%%%%%%%%%%
\part{Background and motivation.}
%%%%%%%%%%%%%%%%%%%%%%%%%%%%%%%%%%
%%%%%%%%%%%%%%%%%%%%%%%%%%%%%%%%%%

We recall the famous integral formula of Gauss \cite{Gauss} for the linking number of two disjoint smooth closed curves
$$X = \{x(s)\,:\,s\in S^1\}\quad\mbox{\rm and}\quad Y = \{y(t)\,:\,t\in S^1\}$$
in Euclidean 3-space $\reals^3$:
$$\lk(X,Y) = \frac{1}{4\pi}\int_{T^2} \frac{dx}{ds} \times \frac{dy}{dt} \dtw \frac{x-y}{|x-y|^3}\, ds\, dt.$$ 
The \textit{\textbf{helicity}} \rm of a vector field $V$ defined on a bounded domain $\Omega$ in $\reals^3$ is
given by the formula 
$$\text{Hel}(V) = \frac{1}{4\pi}\int_{\Omega\times\Omega} V(x)\times V(y) \dtw \frac{x-y}{|x-y|^3}\,dx\,dy\,,$$
where $dx$ and $dy$ are volume elements.

There is no mistaking the analogy with Gauss's linking integral, and no surprise that helicity is a measure of the extent to which the orbits of $V$ wrap and coil around one another.

Woltjer \cite{Woltjer} introduced this notion during his study of the magnetic field in the Crab Nebula, showed that the helicity of a magnetic field remains constant as the field evolves according to the equations of ideal magneto-hydrodynamics, derived from this a lower bound for the changing field energy, and calculated the stable field at the end of the evolution. The term ``helicity'' was coined by Moffatt \cite{Moffatt}, who also derived the above formula from Woltjer's original expression.

Since its introduction, helicity has played an important role in astrophysics and solar physics, and in plasma physics here on earth.

Our study was motivated by a problem proposed by Arnol$'$d and Khesin \cite{ArnoldKhesin98} regarding the search for ``higher helicities'' for divergence-free vector fields. In their own words:

\vspace{.05in}
\begin{indentation}{2.6em}{2.6em}
\small\itb
\noindent The dream is to define such a hierarchy of invariants for generic vector fields such that, whereas all the invariants of order $\leq k$ have zero value for a given field and there exists a nonzero invariant of order $k+1$, this nonzero invariant provides a lower bound for the field energy.
\end{indentation}
\vspace{.05in}

Previous integral formulas for Milnor's triple linking number and attempts to define a higher order helicity can be found in the work of Massey~\cite{Massey58, Massey68}, Monastyrsky and Retakh~\cite{MonastyrskyRetakh}, Berger~\cite{Berger90, Berger91}, Guadagnini, Martellini and Mintchev~\cite{GMM}, Evans andBerger~\cite{EvansBerger92}, Akhmetiev and Ruzmaiken~\cite{Akhmetiev94, Akhmetiev95}, Arnol$'$d and Khesin~\cite{ArnoldKhesin98}, Laurence and Stredulinsky \cite{Laurence}, Leal~\cite{Leal02}, Hornig and Mayer \cite{Hornig}, Rivi\`ere~\cite{Riviere}, Khesin~\cite{Khesin}, Bodecker and Hornig~\cite{Bodecker04}, Auckly and Kapitanski~\cite{AucklyKapitanski}, Akhmetiev~\cite{Akhmetiev05}, and Leal and Pineda~\cite{LealPineda08}.

The principal sources for these formulas are Massey triple products in cohomology, quantum field theory in general, and Chern-Simons theory in particular. A common feature of these integral formulas is that choices must be made to fix the domain of integration and the value of the integrand.

Our own approach to this problem, initiated in Part 1 and continued here, has been applied by Komendarczyk~\cite{Komendarczyk, Komendarczyk10} in special cases to derive a higher order helicity for magnetic fields whose ordinary helicity is zero, and to obtain from this nonzero lower bounds for the field energy.

\clearpage

%%%%%%%%%%%%%%%%%%%%%%%%%%%%%%%%%%
%%%%%%%%%%%%%%%%%%%%%%%%%%%%%%%%%%
\part{The key map from $\conf$ to $S^2$.}
%%%%%%%%%%%%%%%%%%%%%%%%%%%%%%%%%%
%%%%%%%%%%%%%%%%%%%%%%%%%%%%%%%%%%

Let $x$, $y$ and $z$ be three distinct points in Euclidean 3-space $\reals^3$. They will typically span a triangle there, but are permitted to be colinear, as long as they remain distinct.

We show a typical configuration below, with the sides of the triangle oriented and labeled as $a$, $b$, $c$, with the interior angles labeled as $\alpha$, $\beta$, $\gamma$, and with the orientation of the triangle determining a choice of unit normal vector $n$.

If the triangle degenerates to a doubly covered line segment, then the sides are still recognizable, and  likewise the interior angles, with two of them zero and the third $180^\circ$. In this case we set $n$ equal to the zero vector.

%%%% FIGURE 1: Triangle %%%%%
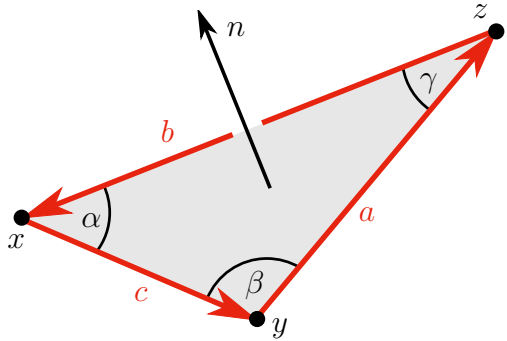
\begin{figure}[h!]
	\includegraphics[height=120pt]{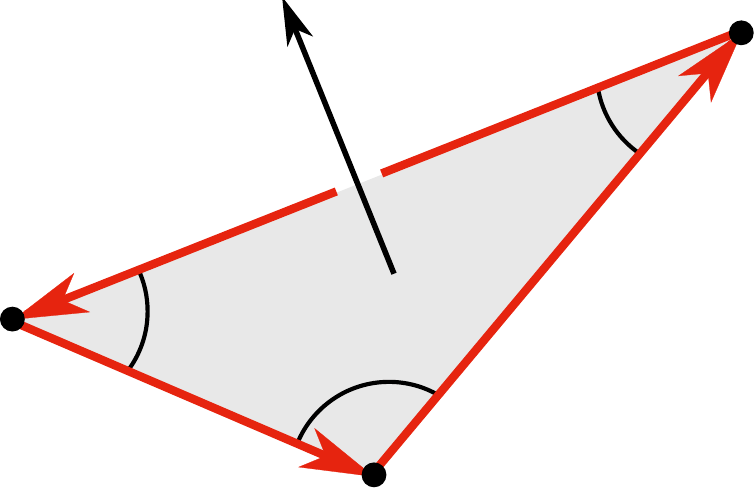}
	\put(-188,30){\large $x$}
	\put(-88,-2){\large $y$}
	\put(-12,117){\large $z$}
	\put(-55,40){\large \textcolor{myred}{$a$}}
	\put(-130,70){\large \textcolor{myred}{$b$}}
	\put(-140,10){\large \textcolor{myred}{$c$}}
	\put(-160,38){\large $\alpha$}
	\put(-98,13){\large $\beta$}
	\put(-32,91){\large $\gamma$}
	\put(-105,110){\large $n$}
	\caption{Triangle and normal vector}
	\label{fig:triangle}
\end{figure}
%%%%%%%%%%%%%%%%%%%

We write $[a] = a /|a|$ for the unit vector pointing along the oriented side $a$ of our triangle, and likewise for the other two sides, and similarly define
$$[b,c] = \frac{b\times c}{|b||c|} = (\sin\alpha)n\,,$$ 
and likewise for the other two pairs of sides.  Next, we define a vector in 3-space by the formula 
\textbf{\boldmath $$F(x,y,z) = [a] + [b] + [c] + [b,c] + [c,a] + [a,b]\,,$$}
equivalently, 
$$F(x,y,z) = \left(\frac{a}{|a|} + \frac{b}{|b|} + \frac{c}{|c|}\right) + (\sin\alpha + \sin\beta + \sin\gamma)n\,.$$

The term $[a] = a /|a|$ is just the classical Gauss map applied to the vertices $y$ and $z$, and the expression $[a] + [b] + [c]$ is the symmetrization of this. It is a vector tangent to the plane of the triangle, and is easily seen to vanish only for equilateral triangles.  

The term $[b,c] = b\times c/|b||c| = (\sin\alpha)n$ is a dimensionless version of the ``directed area'' of the triangle, and the expression $[b, c] + [c, a] + [a, b]$ is the symmetrization of this. This vector is orthogonal to the plane of the triangle, and vanishes only when the triangle degenerates, 
with the vertices lying along a line but remaining distinct.

It follows that $F(x, y, z)$ is never zero, since it is the sum of two orthogonal vectors that do not vanish simultaneously. 

The smoothness of $F$ as a function of $x$, $y$ and $z$ is apparent from its defining formula, since each of its six terms is a smooth function of distinct points. The equivariance of $F$ with respect to orientation preserving isometries of $\reals^3$ is similarly apparent, since
$$[\varphi(a)] = \varphi([a])\quad\mbox{\rm and}\quad [\varphi(b), \varphi(c)] = \varphi([b, c])$$
for any rotation $\varphi$ of $\reals^3$, while translations don't change $a$, $b$ and $c$. $F$ is likewise insensitive to change of scale.

Let $\conf$ denote the \textit{\textbf{configuration space}} of ordered triples of distinct points in $\reals^3$. Then with the above definition, we have
$$F:\conf \to \reals^3 -\{0\}\,.$$ 
Since the image of $F$ misses the origin of $\reals^3$, we may normalize to obtain
$$f = \frac{F}{|F|}:\conf \to S^2\,.$$
Then the map $f$ is also smooth, equivariant as above, and insensitive to change of scale.

%%%%%%%%%%%%%%%%%%%%%%%%%%%%%%%%%%
%%%%%%%%%%%%%%%%%%%%%%%%%%%%%%%%%%
\part{The generalized Gauss map.}
%%%%%%%%%%%%%%%%%%%%%%%%%%%%%%%%%%
%%%%%%%%%%%%%%%%%%%%%%%%%%%%%%%%%%

Suppose now that $L$ is a link in $\reals^3$ with three parametrized components 
$$X = \{x(s)\,:\,s\in S^1\},\quad Y = \{y(t)\,:\,t\in S^1\},\quad Z = \{z(u)\,:\,u\in S^1\}\,.$$
We define the \textit{\textbf{generalized Gauss map}} $g_L : T^3\to S^2$ by 
$$g_L(s, t, u) = f(x(s) , y(t) , z(u))\,.$$
We regard this map as a natural generalization of the classical Gauss map from the 2-torus $T^2$ to the 2-sphere $S^2$ associated with a two-component link in $\reals^3$. If the link $L$ is smooth, then so is the map $g_L$.

The map $g_L$ is equivariant with respect to the group $\Isom^+\reals^3$ of orientation-preserving isometries of $\reals^3$. That is, if $\varphi$ is such an isometry, then $g_{\varphi(L)} = \varphi\circ g_L$, where $\varphi$ acts on $S^2$ via its ``rotational part''. In particular, if $\varphi$ is a translation, then $g_{\varphi(L)} = g_L$. 

Since the map $f$ is insensitive to change of scale, so is the map $g_L$.

The homotopy class of $g_L$ is unchanged under reparametrization of $L$, or more generally under any link homotopy of $L$. The generalized Gauss map is also ``sign symmetric'' in that it transforms under any permutation of the components of $L$ by precomposing with the corresponding permutation automorphism of $T^3$ multiplied by the sign of the permutation.

%%%%%%%%%%%%%%%%%%%%%%%%%%%%%%%%%%
%%%%%%%%%%%%%%%%%%%%%%%%%%%%%%%%%%
\part{Pictures of the generalized Gauss map.}
%%%%%%%%%%%%%%%%%%%%%%%%%%%%%%%%%%
%%%%%%%%%%%%%%%%%%%%%%%%%%%%%%%%%%

In each of the three figures below, we show the vector 
$$F(x,y,z) = \left(\frac{a}{|a|} + \frac{b}{|b|} + \frac{c}{|c|}\right) + (\sin\alpha + \sin\beta + \sin\gamma)n$$
attached to a triangle in $\reals^3$ with vertices at $x$, $y$ and $z$.

\noindent{\itb Equilateral triangle.}

\begin{figure}[h!]
	\subfigure[The vector $F$ for an equilateral triangle.]{\label{fig:eqtriangle}
		\includegraphics[scale=0.9]{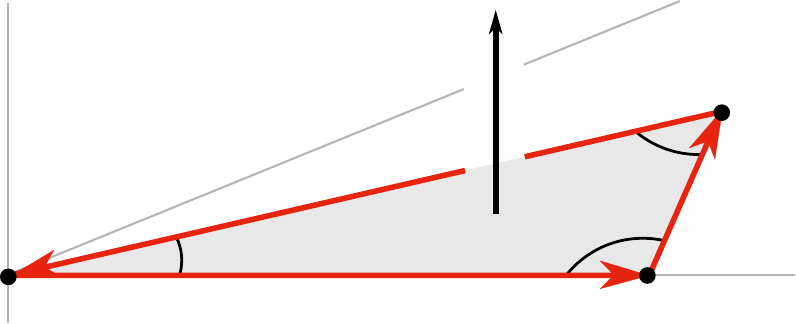}
		\put(-218,11){\large $x$}
		\put(-36,5){\large $y$}
		\put(-15,60){\large $z$}
		\put(-25,30){\large \textcolor{myred}{$a$}}
		\put(-105,42){\large \textcolor{myred}{$b$}}
		\put(-110,5){\large \textcolor{myred}{$c$}}
		\put(-155,17){\tiny $60^\circ$}
		\put(-65,20){\tiny $60^\circ$}
		\put(-50,40){\tiny $60^\circ$}
		\put(-130,70){\large $F \approx 2.6n$}
	}
	\hspace{5pt}
	\subfigure[The vector $F$ for a 30-60-90 right triangle.]{\label{fig:306090triangle}
		\includegraphics[scale=0.9]{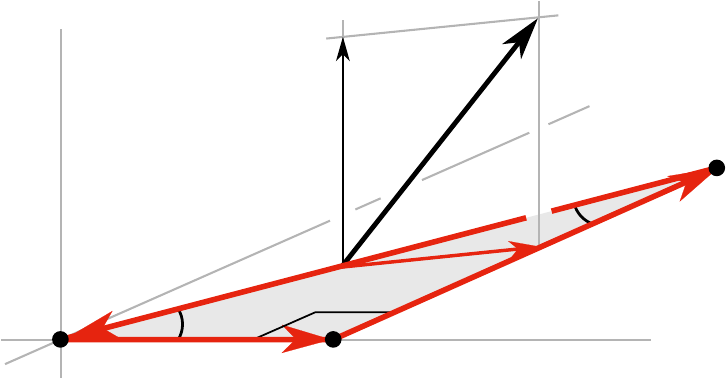}
		\put(-183,13){\large $x$}
		\put(-98,3){\large $y$}
		\put(2,58){\large $z$}
		\put(-70,17){\large \textcolor{myred}{$a$}}
		\put(-118,28){\large \textcolor{myred}{$b$}}
		\put(-138,2){\large \textcolor{myred}{$c$}}
		\put(-137,13){\tiny $60^\circ$}
		\put(-100,20){\tiny $90^\circ$}
		\put(-53,35){\rotatebox{20}{\tiny $30^\circ$}}
		\put(-130,70){\large $2.37n$}
		\put(-47,28){\large \textcolor{myred}{$.52$}}
		\put(-80,70){\large $F$}
	} \\
	\vspace{.3in}
	\subfigure[The vector $F$ for a degenerate triangle.]{\label{fig:degtriangle}
	\includegraphics[scale=0.6]{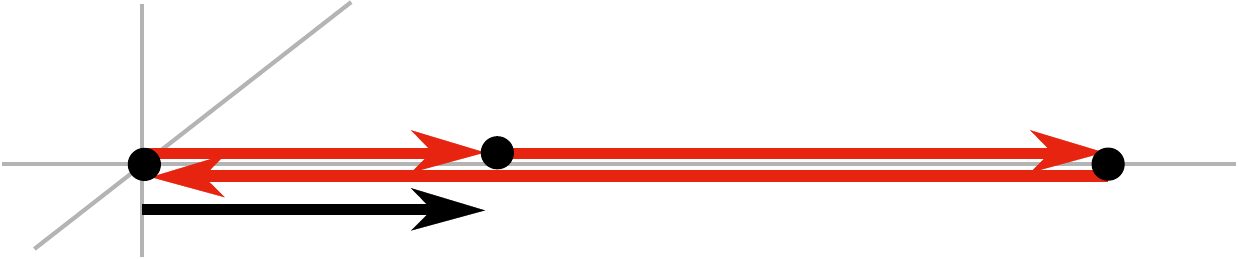}
	\put(-200,20){\large $y$}
	\put(-20,19){\large $x$}
	\put(-132,24){\large $z$}
	\put(-165,23){\large \textcolor{myred}{$a$}}
	\put(-80,23){\large \textcolor{myred}{$b$}}
	\put(-95,5){\large \textcolor{myred}{$c$}}
	\put(-165,-3){\large $F$}
	}
	\label{fig:triangles}
	\caption{The vector $F$ for various triangles}
\end{figure}
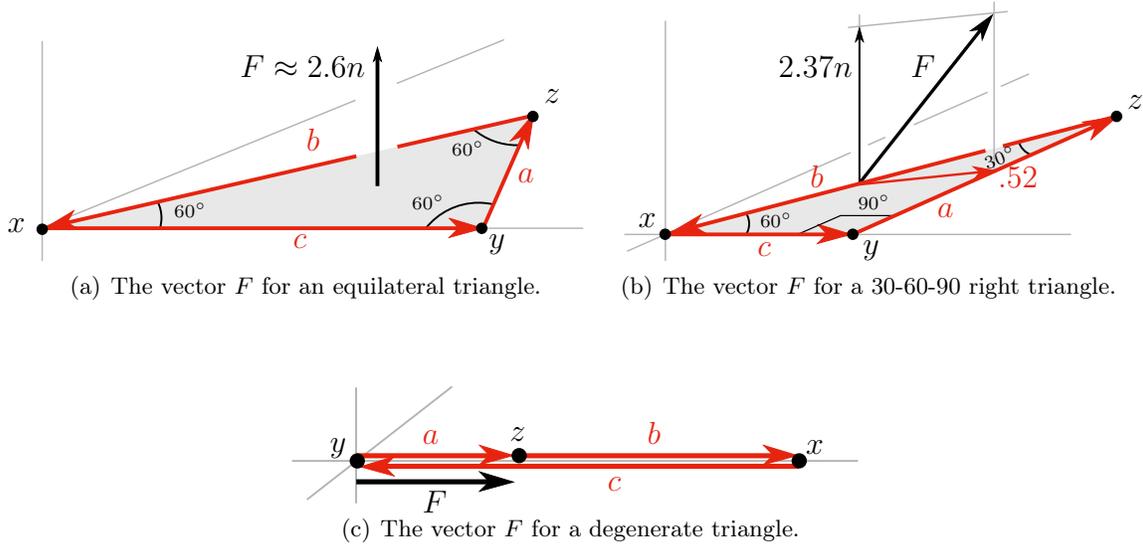

In this case, $(a/|a| + b/|b| + c/|c|) = 0$, and so 
$$F = (\sin60^\circ + \sin60^\circ + \sin60^\circ)n = 3(\sqrt{3}/2)n \approx 2.6\,n,$$
as depicted in Figure~\ref{fig:eqtriangle}.

\clearpage

%%%% FIGURE 3 %%%%%
%%%%%%%%%%%%%%%%%%%

\noindent{\itb 30-60-90 right triangle.}

% %%%% FIGURE 4 %%%%%
% \begin{figure}[h!]
% 	\includegraphics[height=150pt]{306090triangle2.pdf}
% 	\put(-260,6){\large $x$}
% 	\put(-152,6){\large $y$}
% 	\put(0,75){\large $z$}
% 	\put(-95,32){\large \textcolor{myred}{$a$}}
% 	\put(-180,42){\large \textcolor{myred}{$b$}}
% 	\put(-210,8){\large \textcolor{myred}{$c$}}
% 	\put(-227,19){\tiny $60^\circ$}
% 	\put(-163,20){\tiny $90^\circ$}
% 	\put(-56,63){\rotatebox{20}{\tiny $30^\circ$}}
% 	\put(-180,90){\large $2.37n$}
% 	\put(-125,39){\large \textcolor{myred}{$.52$}}
% 	\put(-87,117){\large $F$}
% \end{figure}
% %%%%%%%%%%%%%%%%%%%

In Figure~\ref{fig:306090triangle}, we take $a$ to be of length $\sqrt{3}$, $b$ to be of length $2$, and $c$ to be of length $1$. The vector $(a/|a| + b/|b| + c/|c|)$ is shown above running from a point on the hypotenuse $b$ to a point on the longer side $a$, and has length $\approx .52$. The quantity
$$(\sin90^\circ + \sin60^\circ + \sin30^\circ)n = (1 + \sqrt{3}/2 + 1/2)n \approx 2.37\,n\,,$$ 
and so
$$F \approx \mbox{\rm(horizontal vector of length .52)} + 2.37\,n\,.$$

\noindent{\itb Degenerate triangle.}

% %%%% FIGURE 5 %%%%%
% \begin{figure}[h!]
% 	\includegraphics[scale=0.75]{degtriangle2.pdf}
% 	% \includegraphics[height=60pt]{degtriangle2.pdf}
% 	\put(-249,27){\large $y$}
% 	\put(-23,25){\large $x$}
% 	\put(-165,30){\large $z$}
% 	\put(-200,29){\large \textcolor{myred}{$a$}}
% 	\put(-105,29){\large \textcolor{myred}{$b$}}
% 	\put(-120,8){\large \textcolor{myred}{$c$}}
% 	\put(-205,0){\large $F$}
% 	\caption{ }
% 	\label{fig:degtriangle}
% \end{figure}
% %%%%%%%%%%%%%%%%%%%

The degenerate triangle shown in Figure~\ref{fig:degtriangle} has $a$ of length $1$, $b$ of length $2$, and $c$ of length $3$.  The vector ${(a/|a|+b/|b|+c/|c|)}$ is then of length $1$ lying in the line of the triangle as shown, while
$$(\sin0^\circ + \sin0^\circ + \sin180^\circ)n = 0\,,$$ and so $F$ is the unit vector shown.

%%%%%%%%%%%%%%%%%%%%%%%%%%%%%%%%%%
%%%%%%%%%%%%%%%%%%%%%%%%%%%%%%%%%%
\part{Statement of results.}
%%%%%%%%%%%%%%%%%%%%%%%%%%%%%%%%%%
%%%%%%%%%%%%%%%%%%%%%%%%%%%%%%%%%%

The first of our two main results gives an explicit correspondence between the Milnor link homotopy invariants of a three-component ordered, oriented link in $\reals^3$ and the Pontryagin homotopy invariants of its generalized Gauss map.

%%%%%%%%%%%%%%%%%
\begin{theorem}\label{thm:A}
	\textbf{\boldmath Let $L$ be a three-component link in $\reals^3$. 
Then the pairwise linking numbers $p$, $q$ and $r$ of $L$ are equal to the degrees of 
its generalized Gauss map $g_L\colon T^3\to S^2$ on the two-dimensional coordinate subtori 
of $T^3$, while twice Milnor's $\mu$-invariant for $L$ is equal to 
Pontryagin's $\nu$-invariant for $g_L$ modulo $2\gcd(p, q, r)$.}
\end{theorem}
%%%%%%%%%%%%%%%%%

Each two-dimensional coordinate subtorus of $T^3$ is oriented to have positive intersection with the remaining circle factor.

We refer the reader to Part 1 for a discussion of Milnor's $\mu$-invariant for a three-component link in $\reals^3$ or $S^3$, and of Pontryagin's $\nu$-invariant in the special case of a smooth map of a 3-manifold to the 2-sphere. We also explained there how to convert Pontryagin's original relative invariant to an absolute invariant for maps from the 3-torus to the 2-sphere by comparing with an appropriate family of ``base maps''.

In Part 1, the long and detailed proof of \thmref{A} in the spherical case was carried out in terms of framed bordism of framed links in the 3-torus. We will capitalize on that effort here by building a bridge from the Euclidean to the spherical versions, and cross it to transfer the burden of proof from the Euclidean side back to the spherical side, to work already done there.

Our second main result provides an integral formula for Milnor's $\mu$-invariant in the special case when the pairwise linking numbers $p$, $q$ and $r$ of $L$ vanish.

To explain the symbols that appear in that formula, let $\omega$ denote the usual area form on $S^2$, normalized to have total area $1$. Then $\omega$ pulls back under the generalized Gauss map $g_L\colon T^3\to S^2$ to a closed 2-form $\omega_L$ on $T^3$, which we refer to as the {\itb characteristic 2-form} of $L$.

When $p$, $q$ and $r$ are all zero, it follows from \thmref{A} that $\omega_L$ is exact. We then let $d^{-1}(\omega_L)$ denote any ``primitive'' of $\omega_L$, meaning any 1-form on $T^3$ whose exterior derivative is $\omega_L$.

%%%%%%%%%%%%
\begin{theorem}\label{thm:B}
	\textbf{\boldmath Let $L$ be a three-component link in 
$\reals^3$ whose pairwise linking numbers are all zero. Then Milnor's $\mu$-invariant of $L$ is given by the formula
\begin{equation}\label{eq:B1}
	\mu(L) = \frac12\int_{T^3}	d^{-1}(\omega_L) \wedge \omega_L\, .
\end{equation} 
}
\end{theorem}
%%%%%%%%%%%%

The value of this integral is easily seen to be independent of the choice of
primitive $d^{-1}(\omega_L)$ for $\omega_L$.

The geometrically natural choice for $d^{-1}(\omega_L)$ is the primitive of least $L^2$-norm. It can be obtained explicitly by convolving $\omega_L$ with the fundamental solution $\varphi$ of the scalar Laplacian on $T^3$, and then taking the exterior co-derivative $\delta$ of the resulting 2-form:
$$d^{-1}(\omega_L) = \delta(\varphi \ast \omega_L)\,.$$ 
Details of this construction are presented after the proof of \thmref{B}. If we make this geometrically natural choice for $d^{-1}(\omega_L)$, then the integrand in the above formula is also geometrically natural in the sense that it is unchanged if $L$ is moved by an orientation-preserving isometry of $\reals^3$.

For comparison with formula~\eqref{eq:B1} above, we recall J.H.C.\ Whitehead's integral formula for the Hopf invariant of a smooth map $f\colon S^3\to S^2$,
$$\Hopf(f) = \int_{S^3}d^{-1}(f^\dast\omega) \wedge f^\dast\omega\,,$$
where $\omega$ is the normalized area form on $S^2$, and $f^\dast\omega$ is its pullback via $f$ to an exact 2-form on $S^3$.

\vspace{\baselineskip}

There are two additional versions of the integral formula for Milnor's $\mu$-invariant given in \linebreak \thmref{B}, and we present them next. To state these formulas, we again need some definitions. 

Let $L$ be a three-component link in $\reals^3$, and $\omega_L$ its characteristic 2-form on $T^3$. We convert the closed 2-form $\omega_L$ to a divergence-free vector field $V_L$ on $T^3$ via the usual formula,
$$\omega_L(A,B) = (A\times B)\dtw V_L\,,$$
for all vector fields $A$ and $B$ on $T^3$. We refer to $V_L$ as the \textit{\textbf{characteristic vector field}} of $L$ on $T^3$. When the pairwise linking numbers $p$, $q$ and $r$ of $L$ are all zero, the vector field $V_L$ on $T^3$ is in the image of curl.

For the third version of our formula for Milnor's $\mu$-invariant, we need to express the characteristic 2-form and vector field in terms of Fourier series on the 3-torus. To that end, view $T^3= S^1 \times S^1 \times S^1$ as the quotient $(\reals/2\pi\ints)^3$, and write $\x = (s, t, u)$ for a general point there.

Using the complex form of Fourier series, express 
$$\omega_L = \sum_{\n\in \ints^3} (c_\n^s \,dt\wedge du + 
c_\n^t \,du\wedge ds + c_\n^u \,ds\wedge dt)e^{i\n\dtw\x}\,.$$
We compress notation by writing 
$$\c_\n = (c_\n^s ,c_\n^t ,c_\n^u)\,,$$
$$d\x = (ds,dt,du)\,,\quad \star d\x = (dt\wedge du\,,\,du\wedge ds\,,\,ds\wedge dt)\,,$$
$$\partial_\x = (\partial_s ,\partial_t ,\partial_u) = \left(\frac{\partial}{\partial s},
\frac{\partial}{\partial t},\frac{\partial}{\partial u}\right)\,.$$
Using this compression, the formulas for $\omega_L$ and $V_L$ become 
$$\omega_L = \sum_\n \c_\n e^{i\n\dtw\x}\dtw¥\star dx\quad
\mbox{\rm and\ } V_L = \sum_\n \c_\n e^{i\n\dtw\x}\dtw \partial_\x\,.$$
Writing $\zero = (0,0,0)$, the coefficient $\c_\zero = (c_\zero^s ,c_\zero^t ,c_\zero^u) = \zero$ since the form $\omega_L$ is exact, equivalently the vector field $V_L$ is in the image of curl. Finally, we express the general Fourier coefficient $\c_\n$ in terms of its real and imaginary parts,
$$\c_\n =\a_\n +i\b_\n\, \quad\mbox{\rm with $\a_\n$ and $\b_\n$ real.}$$

\begin{thmBcont}
	\textbf{\boldmath Let $L$ be a three-component link in $\reals^3$ with pairwise linking numbers all zero. Then Milnor's $\mu$-invariant of $L$ is also given by the formulas
\begin{align}
	\mu(L) \ &=\ \frac12\int_{T^3\times T^3} V_L(x)\times V_L(y)\dtw\nabla_y\varphi(x-y)\,dx\,dy   \tag{2} \label{eq:B2}\\
\ &  = \  8\pi^3\sum_{\n\neq\zero}\a_n\times\b_\n\dtw\frac{\n}{\Vert\n\Vert^2}\,, \tag{3} \label{eq:B3}
\end{align}
where $\varphi$ is the fundamental solution of the scalar Laplacian on the 3-torus, $V_L$ is 
the characteristic vector field of $L$, and $\a_\n$ and $b_\n$ are the real and imaginary parts of 
the Fourier coefficients $\c_\n$ of $V_L$.
}
\end{thmBcont}

In formula \eqref{eq:B2} above, the difference $x-y$ is taken in the abelian group structure on $T^3$, the expression $\nabla_y\varphi(x-y)$ indicates the gradient with
respect to $y$ while $x$ is held fixed, and $dx$ and $dy$ are volume elements on $T^3$.

Formula \eqref{eq:B2} is just the vector field version of formula \eqref{eq:B1}, in which the integral hidden in the convolution formula for $d^{-1}(\omega_L)$ is expressed openly. This formula shows that the Milnor triple linking number $\mu(L)$ is one-half the helicity of the vector field $V_L$ on the 3-torus $T^3$. The integrand in formula \eqref{eq:B2} is invariant under the group $\Isom^+\reals^3$ of orientation-preserving isometries of $\reals^3$.

\clearpage

%%%%%%%%%%%%%%%%%%%%%%%%%%%%%%%%%%
%%%%%%%%%%%%%%%%%%%%%%%%%%%%%%%%%%
\part{Numerical computation.}
%%%%%%%%%%%%%%%%%%%%%%%%%%%%%%%%%%
%%%%%%%%%%%%%%%%%%%%%%%%%%%%%%%%%%
% 
% \medskip
% 
% \noindent\bf Numerical computation \rm  

\begin{figure}[h!]
\includegraphics[height=170pt]{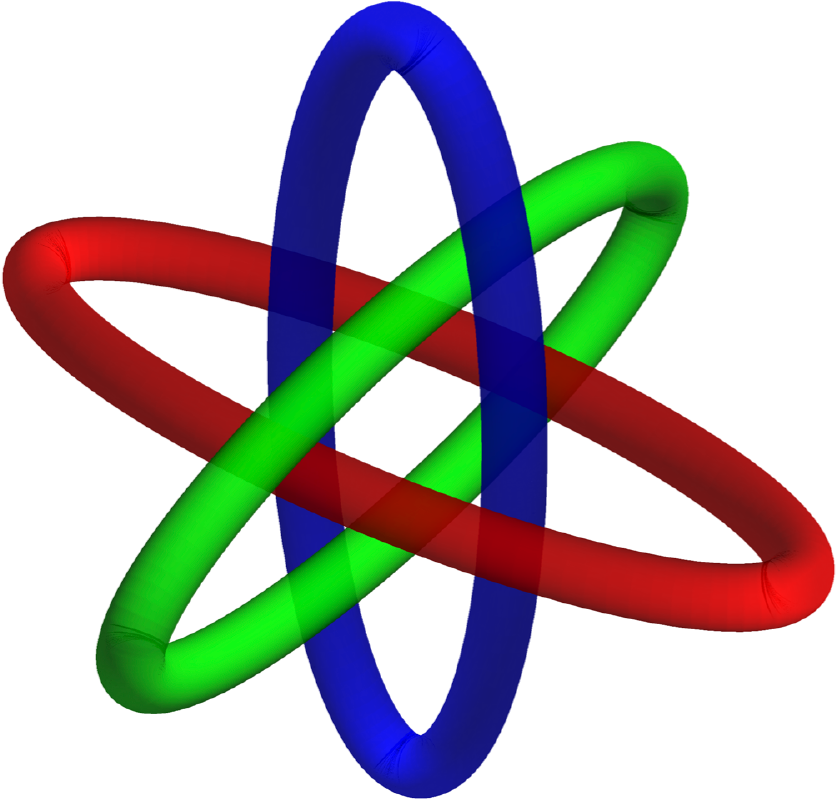}
\caption{Borromean rings}
\label{fig:borrcalc}
\end{figure}

We used Matlab to calculate an approximation to Milnor's $\mu$-invariant, as given by formula \eqref{eq:B3} of \thmref{B}, for the three-component link $L$ in $\reals^3$ parametrized by 
\begin{align*}
x(s) &= \textstyle{(2\cos s, 7\sin s, 0)}\cr
y(t) &= \textstyle{(0,2\cos t, 7\sin t)}\cr
z(u) &= \textstyle{(7\sin u, 0, 2\cos u)}\,,
\end{align*}
with $s\in [0,2\pi]$, $t\in [0,2\pi]$, $u\in[0,2\pi]$, which is a concrete realization of the Borromean rings with $\mu = -1$ shown in Figure~\ref{fig:borrcalc}.

%%%% FIGURE 6 %%%%%
%%%%%%%%%%%%%%%%%%%

In particular, we used Matlab to calculate approximations to the Fourier coefficients $\c_\n$ of its characteristic form $\omega_L$. We used subdivisions of the $s$, $t$ and $u$ intervals into 256 subintervals to approximate the integrals defining the coefficients $\c_\n$ for $-64 \le n_s,n_t,n_u\le 64$.  The approximation of $\mu$ we obtained in this way was $-0.99999997$.

% \clearpage
%%%%%%%% THIS IS THE END OF THE INTRODUCTION %%%%%%%%

%% file: 2thmA.tex
%%%%%%%%%%%%%%%%%%%%%%%%%
%% SECTION 2: Theorem A
%%%%%%%%%%%%%%%%%%%%%%%%%

\section{\thmref{A}} 
\label{sec:thmA}

%%%%%%%%%%%%%%%%%%%%%%%%%%%%%
%%%%%%%%%%%%%%%%%%%%%%%%%%%%%
\part{Proof plan for \thmref{A}.}
%%%%%%%%%%%%%%%%%%%%%%%%%%%%%
%%%%%%%%%%%%%%%%%%%%%%%%%%%%%

We begin with a quick summary of the spherical theory, and its key map $f_S\colon \confs \to S^2$. Next, we discuss inverse stereographic projection $h$ from $\reals^3$ to $S^3$, and use it to define a map $H\colon \conf \to \confs$.  Afterwords, we state and prove the ``bridge lemma'', which asserts the homotopy commutativity of the diagram
\vskip-.1in
\begin{diagram}[size=2em]
	\confs & \rTo[l>=.5in]^{f_S} & S^2 \\
	\uTo^{H} & & \dEquals \\
	\conf & \rTo^{f_E} & S^2
\end{diagram}
\vskip.1in
\noindent where the horizontal maps are the ``key maps'' of the spherical and Euclidean theories. Finally, we use the bridge lemma to prove \thmref{A}.

%%%%%%%%%%%%%%%%%%%%%%%%%%%%%
%%%%%%%%%%%%%%%%%%%%%%%%%%%%%
\part{Recollection of the spherical theory.}
%%%%%%%%%%%%%%%%%%%%%%%%%%%%%
%%%%%%%%%%%%%%%%%%%%%%%%%%%%%

The key map $\confs \to S^2$ in the spherical setting was defined in Part 1 as follows. Let $x$, $y$ and $z$ be three distinct points on the unit 3-sphere $S^3$ in $\reals^4$. They cannot lie on a straight line in $\reals^4$, so must span a 2-plane there. Translate this plane to pass through the origin, and then orient it so that the vectors $x-z$ and $y-z$ form a positive basis. The result is an element $G(x, y, z)$ of the Grassmann manifold $G_2\reals^4$ of all oriented 2-planes through the origin in 4-space. This procedure defines the {\itb Grassmann map}
$$G\colon \confs\to G_2\reals^4\,,$$ 
where $\confs$ is the configuration space of ordered triples of distinct points in $S^3$ .

%%%% FIGURE 7: Grassmann Map %%%%%
\begin{figure}[h!]
	\includegraphics[height=120pt]{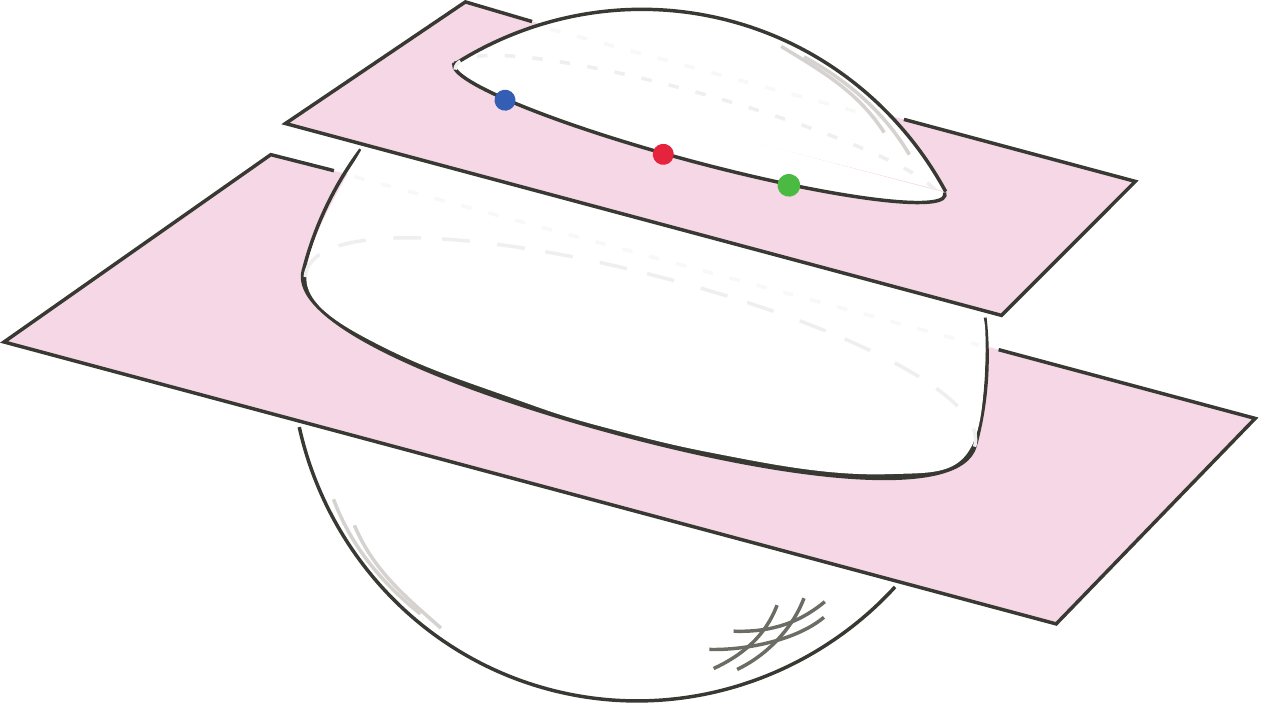}
	\put(-70,27){\small $G(x,y,z)$}
	\put(-135,95){\small $x$}
	\put(-107,86){\small $y$}
	\put(-85,80){\small $z$}
	\put(-125,15){$S^3$}
	\put(-190,23){$\br^4$}
\caption{{Map from $\confs$ to $G_2\reals^4$}}
\end{figure}
%%%%%%%%%%%%%%%%%%%

The Grassmann manifold $G_2\reals^4$ with its natural Riemannian metric is, up to scale, isometric to the product $S^2\times S^2$ of two unit 2-spheres. We will express this by the map $\pi\colon G_2\reals^4 \to S^2\times S^2$ which takes the oriented 2-plane $\langle a , b\rangle$ with orthonormal basis $a,b$ to the point $(b \bah{a} , \bah{a} b)$ in $S^2\times S^2$, using quaternion notation and conjugation. This gives us two projection maps $\pi_+$ and $\pi_-$ from $G_2\reals^4\to S^2$,
$$\pi_+\langle a,b\rangle =b\bah{a}\quad\mbox{\rm and}\quad\pi_-\langle a,b\rangle = \bah{a}b\,.$$ 
If the basis $a$, $b$ is not necessarily orthonormal, then we saw in Part 1 that 
$$\pi_+\langle a , b \rangle = \frac{\Im(b\bah{a})}{\Vert \Im(b\bah{a})\Vert}
\quad\mbox{\rm and}\quad \pi_-\langle a , b \rangle =\frac{\Im(\bah{a}b)}{\Vert\Im(\bah{a} b)\Vert}\,.$$ 
We arbitrarily use the first projection $\pi_+$ to define the key map
$$f_S = \pi_+\circ G\colon\confs\to S^2\,.$$

%%%%%%%%%%%%%%%%%%%%%%%%%%%%%
%%%%%%%%%%%%%%%%%%%%%%%%%%%%%
\part{Inverse stereographic projection.}
%%%%%%%%%%%%%%%%%%%%%%%%%%%%%
%%%%%%%%%%%%%%%%%%%%%%%%%%%%%

The corresponding key map in the Euclidean theory is
$$f_E = \frac{F_E}{|F_E|}\colon \conf\to S^2\,,$$
where 
$$F_E(x,y,z) = \left(\frac{a}{|a|} + \frac{b}{|b|} + \frac{c}{|c|}\right) 
+ (\sin\alpha + \sin\beta + \sin\gamma)n$$
was defined earlier, and where we have added the subscripts to 
$f_E$ and $F_E$ to signal ``Euclidean''.

\clearpage

The bridge between the two theories will be a map $H\colon \conf\to\confs$ which makes the diagram
\vskip-.1in
\begin{diagram}[size=2em]
	\confs & \rTo[l>=.5in]^{f_S} & S^2 \\
	\uTo^{H} & & \dEquals \\
	\conf & \rTo^{f_E} & S^2
\end{diagram}
\vskip.1in
\noindent commutative up to homotopy. We present two versions of $H$, the first straightforward via inverse stereographic projection, and the second homotopic to it but more convenient for our arguments.

Viewing $\reals^4$ as the space of quaternions, we regard $\reals^3$ as the subspace of purely imaginary quaternions, and then use inverse stereographic projection from $-1$ to provide a diffeomorphism $h\colon \reals^3\to S^3-\{-1\}$, which preserves the usual orientations on $\reals^3$ and $S^3$.

%%%% FIGURE 8 %%%%%
\begin{figure}[h!]
	\includegraphics[height=165pt]{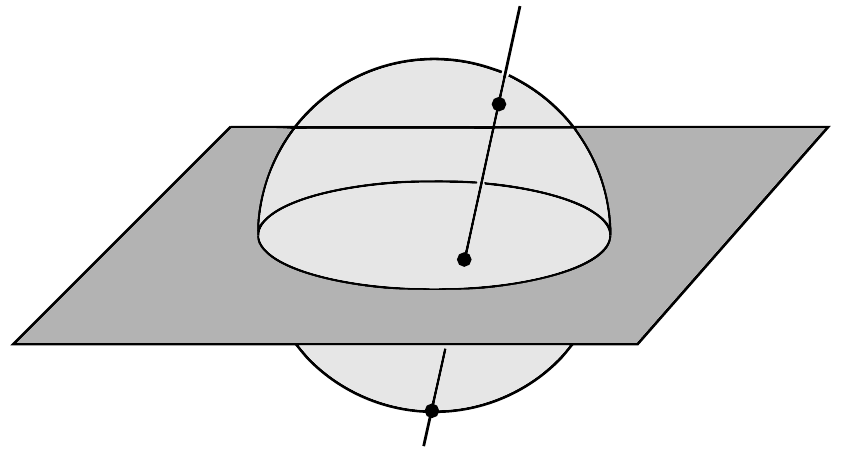}
	\put(-270,47){\Large $\br^3$}
	\put(-185,142){\Large $S^3$}
	\put(-145,7){$-1$}
	\put(-130,70){$q$}
	\put(-146,125){$h(q)$}
	\put(-67,104){$ijk$-space}
	\caption{Inverse stereographic projection $h\colon\reals^3\longrightarrow S^3-\{-1\}$}
\end{figure}
%%%%%%%%%%%%%%%%%%%

Let $q$ denote a purely imaginary quaternion, thus a point of $\reals^3$. We compute that
$$h(q)=\frac{1-|q|^2}{1+|q|^2} + \frac{2q}{1+|q|^2}\,,$$
with the first term on the right being the real part of $h(q)$, and the second term its imaginary part. Indeed, a quick check shows that $h(q)$ has norm 1, and that $h(q)- (-1)$ is a real multiple of $q - (-1)$, and hence that the points $-1$ , $q$ and $h(q)$ lie on a straight line.

The first version of the map $H\colon \conf\to\confs$ uses inverse stereographic projection on each of three points,
$$H(x, y, z) = (h(x), h(y), h(z))\,.$$
The second version, call it $H'$, is defined as follows. Let $(x, y, z)$ be a triple of distinct points in $\reals^3$, and use translation by $-z$ there to move this to the triple $(x-z,y-z,0)$ of distinct points. Then apply $h$ to each of the three points in this new triple to obtain 
$$H'(x,y,z) = (h(x-z),h(y-z),1)\,.$$
The maps $H$ and $H'\colon \conf\to\confs$ are clearly homotopic.

\clearpage

%%%% FIGURE 9 %%%%%
\begin{figure}[h!]
	\includegraphics[height=150pt]{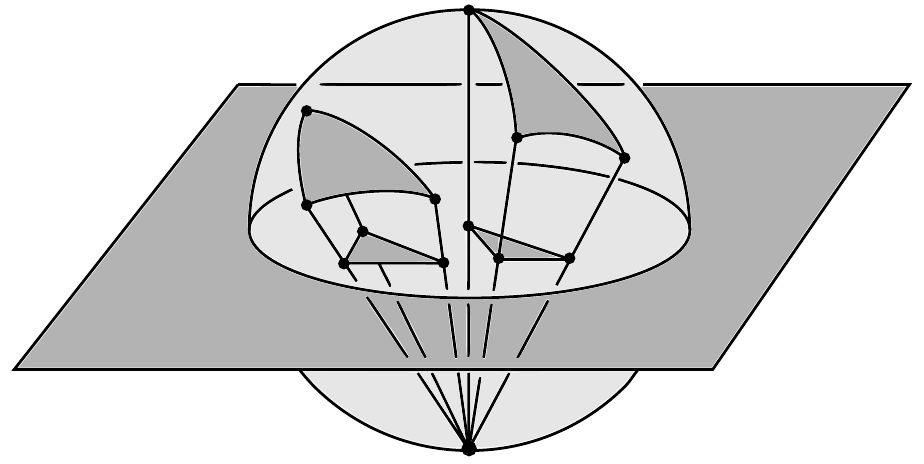}
	\put(-195,140){\Large $S^3$}
	\put(-45,105){\Large $\br^3$}
	\put(-150,150){\tiny $1$}
	\put(-194,63){\tiny $x$}
	\put(-153,64){\tiny $y$}
	\put(-188,74){\tiny $z$}
	\put(-215,77){\tiny $h(x)$}
	\put(-175,80){\tiny $h(y)$}
	\put(-196,116){\tiny $h(z)$}
	\put(-153,75){\tiny $0$}
	\put(-137,61){\tiny $x-z$}
	\put(-111,66){\tiny $y-z$}
	\put(-132,98){\tiny $h(x-z)$}
	\put(-101,113){\rotatebox{-61}{\tiny $h(y-z)$}}
	\put(-192,93){\rotatebox{18}{$H$}}
	\put(-130,112){\rotatebox{27}{$H'$}}
	\caption{The maps $H$ and $H'\colon \conf\longrightarrow\confs$}
\end{figure}
%%%%%%%%%%%%%%%%%%%

%%%%%%%%%%%%%%%%%%%%%%%%%%%%%
%%%%%%%%%%%%%%%%%%%%%%%%%%%%%
\part{Statement of the bridge lemma.}
%%%%%%%%%%%%%%%%%%%%%%%%%%%%%
%%%%%%%%%%%%%%%%%%%%%%%%%%%%%

The result below will permit us to transfer the burden of proof for 
our current Euclidean version of \thmref{A} back to its spherical version in Part 1.

\begin{bridge}
	\textbf{\boldmath The map $H\colon \conf\to\confs$ makes the diagram
\vskip-.1in
\begin{diagram}[size=2em]
	\confs & \rTo[l>=.5in]^{f_S} & S^2 \\
	\uTo^{H} & & \dEquals \\
	\conf & \rTo^{f_E} & S^2
\end{diagram}
\vskip.1in
\noindent commutative up to homotopy.}
\end{bridge}

{\itb Setup.} We avoid the nuisance of normalization by using instead 
the maps 
$$F_S\colon\confs\to\reals^3-\{0\}\quad\mbox{\rm and}\quad
F_E\colon\conf\to\reals^3-\{0\}$$
defined by 
$$F_S(u,v,w)=\Im\!\left((v-w)\overline{(u-w)}\right)$$
and
$$F_E(x,y,z) = \left(\frac{a}{|a|} + \frac{b}{|b|} + \frac{c}{|c|}\right) + 
(\sin\alpha + \sin\beta + \sin\gamma)n\,.$$ 
At the same time, we replace the map $H$ by the homotopic map $H'$. 
So now our job is to show homotopy commutativity of the diagram
\vskip-.1in
\begin{diagram}[size=2em]
	\confs & \rTo[l>=.5in]^{F_S} & \br^3 - \{0\} \\
	\uTo^{H'} & & \dEquals \\
	\conf & \rTo^{F_E} & \br^3 - \{0\}
\end{diagram}
\vskip.1in

\clearpage

%%%%%%%%%%%%%%%%%%%%%%%%%%%%%
%%%%%%%%%%%%%%%%%%%%%%%%%%%%%
\part{Proof of the Bridge Lemma.}
%%%%%%%%%%%%%%%%%%%%%%%%%%%%%
%%%%%%%%%%%%%%%%%%%%%%%%%%%%%

We start with three distinct points $x$, $y$ and $z$ in $\reals^3$, forming a possibly degenerate triangle with sides $a = z-y$, $b = x-z$ and $c = y-x$. Then we begin to compute, 
\begin{align*}
F_S \circ H'(x, y, z) & = F_S(h(x-z) , h(y-z) , 1)\\
& = F_S(h(b) , h(-a) , 1) \\
& = \Im\left((h(-a)-1)\overline{(h(b)-1)}\right)\,.
\end{align*}
We recall the formula
$$h(q)=\frac{1-|q|^2}{1+|q|^2} + \frac{2q}{1+|q|^2}\,,$$
and first substitute $-a$ for $q$, and then $b$ for $q$, to get
$$h(-a)-1=-\frac{2(|a|^2+a)}{1+|a|^2}$$
and
$$\overline{h(b)-1}=-\frac{2(|b|^2+b)}{1+|b|^2}\,.$$
Then
$$\Im\!\left((h(-a)-1)\overline{(h(b)-1)}\right)=C\Im\!\left((|a|^2+a)(|b|^2+b)\right)\,$$
where the positive real number $C$ is given by
$$C = 4\,(1+|a|^2)^{-1}(1+|b|^2)^{-1}\,.$$
We keep in mind that $a$ and $b$, since they lie in $\reals^3$, are purely imaginary quaternions, and hence the sums $|a|^2 + a$ and $|b|^2 + b$ are both quaternions written in terms of their real and imaginary parts.

We recall the following formula about quaternion multiplication, 
$$\Im(q_1 q_2) = \Re(q_1) \Im(q_2) + \Im(q_1) \Re(q_2) + \Im(q_1)\times \Im(q_2)\,,$$
using the vector cross product in the 3-space of purely imaginary quaternions. 

Applying this formula, we get
$$\Im\!\left((|a|^2+a)(|b|^2+b)\right)= |a|^2b+a|b|^2+a\times b\,.$$
Then, stringing together the above computations, we have shown that
$$F_S\circ H'(x,y,z)=C(|a|^2b+a|b|^2+a\times b)\,$$
with
$$C = 4\,(1+|a|^2)^{-1}(1+|b|^2)^{-1}\,.$$
Since we are focusing on homotopy of maps into $\reals^3-\{0\}$,
the strictly positive quantity $C$ is irrelevant, and we hide it from view by recording the 
homotopy
\begin{equation}\label{eq:Fhomotopy}\tag{1}
	F_S\circ H'(x,y,z)\simeq |a|^2b+a|b|^2+a\times b\,.
\end{equation}

It remains to show that this expression is homotopic in $\reals^3-\{0\}$ to
\begin{equation}\label{eq:FE}\tag{2}
	F_E(x,y,z) = \left(\frac{a}{|a|} + \frac{b}{|b|} + \frac{c}{|c|}\right) + 
(\sin\alpha + \sin\beta + \sin\gamma)n\,.
\end{equation}

The right side of \eqref{eq:Fhomotopy} can be rewritten as
\begin{equation}\label{eq:Fhomotopy2}\tag{$1'$}
	\left(|a|^2b+a|b|^2\right)+\left(|a||b|\sin\gamma\right)n
\end{equation}
for convenience of comparison with \eqref{eq:FE}. In each case we have the sum of a vector parallel to the plane of the triangle $xyz$ and a vector orthogonal to it. 

As long as the triangle is non-degenerate, the components in \eqref{eq:Fhomotopy2} and \eqref{eq:FE} orthogonal to its plane are both strictly positive multiples of the unit normal vector $n$. Hence \eqref{eq:Fhomotopy2} and \eqref{eq:FE} are vectors which both lie in the same open half-space of $\reals^3$, and so the line segment between them misses the origin. Thus $F_S\circ H'$ is homotopic to $F_E$ in such cases.

So the issue now is, what happens when the triangle degenerates? In such a case, \eqref{eq:Fhomotopy2} and \eqref{eq:FE} reduce to their tangential components, 
\begin{equation}\tag{$1'$}
	|a|^2b+a|b|^2
\end{equation}
and
\begin{equation}\tag{2}
	\frac{a}{|a|} + \frac{b}{|b|} + \frac{c}{|c|}\,.
\end{equation}
Both of these vectors point along the line of the degenerate triangle, and we must check that they always point the same way, so that the line segment between them again misses the origin.

%%%% FIGURE 10 %%%%%
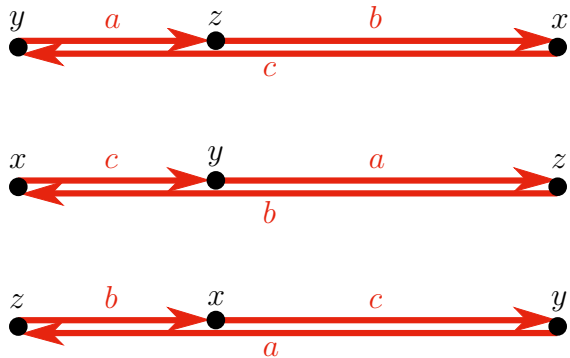
\begin{figure}[h!]
	\includegraphics[height=120pt]{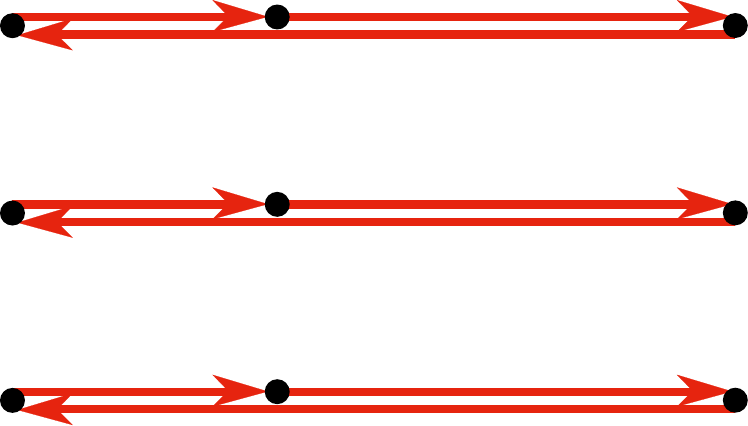}
	\put(-211,121){\large $y$}
	\put(-136,121){\large $z$}
	\put(-6,120){\large $x$}
	\put(-175,120){\large \textcolor{myred}{$a$}}
	\put(-75,120){\large \textcolor{myred}{$b$}}
	\put(-115,102){\large \textcolor{myred}{$c$}}
	\put(-211,67){\large $x$}
	\put(-136,70){\large $y$}
	\put(-6,67){\large $z$}
	\put(-175,67){\large \textcolor{myred}{$c$}}
	\put(-75,67){\large \textcolor{myred}{$a$}}
	\put(-115,46){\large \textcolor{myred}{$b$}}
	\put(-211,14){\large $z$}
	\put(-136,15){\large $x$}
	\put(-6,15){\large $y$}
	\put(-175,14){\large \textcolor{myred}{$b$}}
	\put(-75,14){\large \textcolor{myred}{$c$}}
	\put(-115,-4){\large \textcolor{myred}{$a$}}
	\caption{Degenerate triangles}
	\label{fig:detriangles}
\end{figure}
%%%%%%%%%%%%%%%%%%%

There are three cases, according as which vertex is between the other two. They are shown in Figure~\ref{fig:detriangles}, in which the line of the degenerate triangle is turned so that the two shorter sides point to the right.

In each case, the vector $a/|a| + b/|b| + c/|c|$ is a unit vector pointing to the right. Furthermore, in all three cases, the vector $|a|^2 b + a |b|^2$ points in the same direction as its rescaling $b / |b|^2 + a / |a|^2$, which is a nonzero vector pointing in the same direction as the shorter of the two vectors $a$ and $b$, and this is also to the right. It follows that in every case, non-degenerate or degenerate, the line segment between the vectors \eqref{eq:Fhomotopy2} and \eqref{eq:FE} misses the origin.

Hence the maps $F_S\circ H'$ and $F_E$ from $\conf\to\reals^3-\{0\}$
are homotopic, completing the proof of the Bridge Lemma.

%%%%%%%%%%%%%%%%%%%%%%%%%%%%%
%%%%%%%%%%%%%%%%%%%%%%%%%%%%%
\part{Proof of \thmref{A}.}
%%%%%%%%%%%%%%%%%%%%%%%%%%%%%
%%%%%%%%%%%%%%%%%%%%%%%%%%%%%

Let $L$ be a three-component link in $\reals^3$, let $h\colon\reals^3\to S^3-\{-1\}$ be inverse-stereographic projection, and let $h(L)$ be the resulting three-component link in $S^3$. Since $h$ is an orientation-preserving diffeomorphism, the Milnor invariants $p$, $q$, $r$ and $\mu$ for $L$ match those for $h(L)$.

The Euclidean generalized Gauss map $g_L\colon T^3\to S^2$ for $L$ is given by 
$$g_L(s, t, u) = f_E\left(x(s) , y(t) , z(u)\right)\,,$$
while the spherical generalized Gauss map 
$g_{h(L)}\colon T^3\to S^2$ for $h(L)$ is given by
\begin{align*}
	g_{h(L)}(s, t, u) & = f_S\left(h(x(s)) , h(y(t)) , h(z(u))\right)\\
	& = f_S\circ H\left(x(s),y(t),z(u)\right)\,.
\end{align*}
According to the Bridge Lemma, the maps $f_E$ and $f_S\circ H\colon \conf\to S^2$ are homotopic. 

It follows that the maps $g_L$ and $g_{h(L)}\colon T^3\to S^2$ are also homotopic. Hence the Pontryagin invariants $p$, $q$, $r$ and $\nu$ for $g_L$ match those for $g_{h(L)}$.

Then the correspondence between the Milnor invariants for $L$ and the Pontryagin 
invariants for $g_L$, as asserted in our current Euclidean version of \thmref{A}, 
follows from the correspondence between these invariants for $h(L)$ and 
$g_{h(L)}$, as asserted in the spherical version of \thmref{A}, which was 
proved in Part 1.

This completes the proof of the Euclidean version of \thmref{A}.

%%%%%%%% THIS IS THE END OF SECTION 2 %%%%%%%%

%% file: 3thmB.tex
%%%%%%%%%%%%%%%%%%%%%%%%%
%% SECTION 3: THEOREM B
%%%%%%%%%%%%%%%%%%%%%%%%%

\section{\thmref{B}} 
\label{sec:thmB}

%%%%%%%%%%%%%%%%%%%%%%%%%%%%%
%%%%%%%%%%%%%%%%%%%%%%%%%%%%%
\part{Proof plan for \thmref{B}.}
%%%%%%%%%%%%%%%%%%%%%%%%%%%%%
%%%%%%%%%%%%%%%%%%%%%%%%%%%%%

Let $L$ be a three-component link in Euclidean space $\reals^3$, with pairwise linking numbers all zero. \thmref{B} gives three explicit formulas for the triple linking number (Milnor $\mu$-invariant) of $L$:
\begin{align}
	\tag{1} \mu(L) & = \frac12\int_{T^3}	d^{-1}(\omega_L) \wedge \omega_L \\
	\tag{2}  & = \frac12\int_{T^3\times T^3} V_L(x)\times V_L(y)\dtw\nabla_y\varphi(x-y)\,dx\,dy \\
 	\tag{3} & = 8\pi^3\sum_{\n\neq\zero}\a_n\times\b_\n\dtw\frac{\n}{\Vert\n\Vert^2}\,,
\end{align}
using the notation defined after the two statements of \thmref{B}.

We saw in Part 1 that the spherical version of \thmref{A} implies the spherical version of \thmref{B}, and we show below that the same implication holds for the Euclidean versions here. We will give details only for the proof of formula \eqref{eq:B1}, and refer the reader to Part 1 for the derivation of formulas \eqref{eq:B2} and \eqref{eq:B3}.

The first step will be to give an explicit formula for the characteristic 2-form 
$\omega_L$ of the link $L$.

As mentioned earlier, the geometrically natural choice for $d^{-1}(\omega_L)$ is the primitive of least $L^2$-norm, which can be obtained explicitly by convolving $\omega_L$ with the fundamental solution $\varphi$ of the scalar Laplacian on $T^3$, and then taking the exterior co-derivative $\delta$ of the resulting 2-form: $d^{-1}(\omega_L) = \delta (\varphi\ast\omega_L)$. 

So we will give the expression for this fundamental solution $\varphi$, and then show how formula \eqref{eq:B1} follows from J.H.C.\ Whitehead's integral formula for the Hopf invariant.

%%%%%%%%%%%%%%%%%%%%%%%%%%%%%
%%%%%%%%%%%%%%%%%%%%%%%%%%%%%
\part{An explicit formula for the characteristic 2-form $\omega_L$.}
%%%%%%%%%%%%%%%%%%%%%%%%%%%%%
%%%%%%%%%%%%%%%%%%%%%%%%%%%%%

We start with a three-component link $L$ in $\reals^3$ with components
$$X = \{x(s)\,:\,s\in S^1\},\quad Y = \{y(t)\,:\,t\in S^1\},\quad Z = \{z(u)\,:\,u\in S^1\}\,,$$ 
and recall the figure

%%%% FIGURE 11: Triangle %%%%%
\begin{figure}[h!]
	\includegraphics[height=120pt]{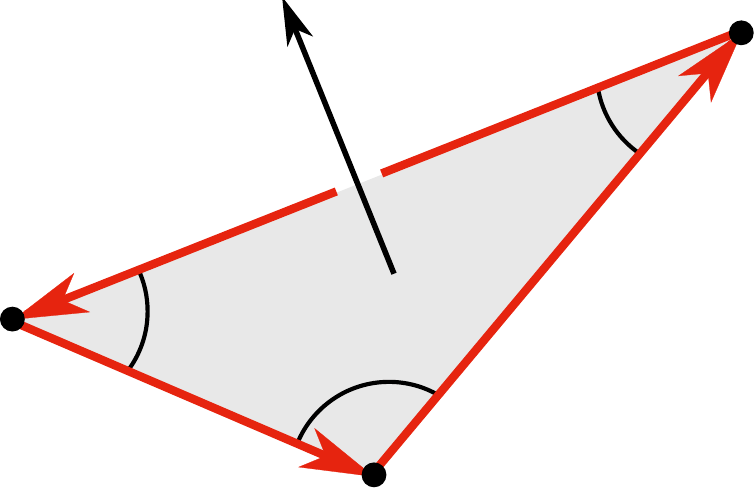}
	\put(-188,30){\large $x$}
	\put(-88,-2){\large $y$}
	\put(-12,117){\large $z$}
	\put(-55,40){\large \textcolor{myred}{$a$}}
	\put(-130,70){\large \textcolor{myred}{$b$}}
	\put(-140,10){\large \textcolor{myred}{$c$}}
	\put(-160,38){\large $\alpha$}
	\put(-98,13){\large $\beta$}
	\put(-32,91){\large $\gamma$}
	\put(-105,110){\large $n$}	% 
		% \caption{Triangle and normal vector}
		% \label{fig:triangle}
\end{figure}
%%%%%%%%%%%%%%%%%%%
\noindent with 
$$a(t,u) = z(u)-y(t),\quad b(s,u) = x(s)-z(u),\quad c(s,t) = y(t)-x(s)\,,$$
and the formula 
$$F(x,y,z) = [a] + [b] + [c] + [b,c] + [c,a] + [a,b]\,.$$
With mild abuse of notation, we write 
$$F(s, t, u) = F(x(s), y(t), z(u))\,.$$
The generalized Gauss map $g_L\colon T^3\to S^2$ of the link $L$ is then given by 
$$g_L(s,t,u) = \frac{F(s,t,u)}{|F(s,t,u)|}\,.$$

Let $\omega$ be the Euclidean area 2-form on the unit 2-sphere $S^2\subset\reals^3$, normalized so that the total area is 1 instead of $4\pi$. If $P$ is a point of $S^2$, and 
$A$ and $B$ are tangent vectors to $S^2$ at $P$, then
$$\omega_P(A,B)=\frac{1}{4\pi}(A\times B)\dtw P\,.$$
This 2-form $\omega$ on $S^2$ extends to a closed 2-form $\overline{\omega}$ on
$\reals^3-\{0\}$ given by
$$\overline{\omega}_P(A,B) = \frac{(A\times B)\dtw P}{4\pi|P|^3}\,,$$ 
which is the pullback of $\omega$ from $S^2$ to $\reals^3-\{0\}$ via the map $P\mapsto P/|P|$.

Hence the pullback $g_L^\dast\omega$ of $\omega$ from $S^2$ to $T^3$ via 
$g_L = F/|F|$ is the same as the pullback $F^\dast\overline{\omega}$ of $\overline{\omega}$ from $\reals^3-\{0\}$ to $T^3$ via $F$. Write 
$$g_L^\dast\omega= F^\dast\overline{\omega}= 
p(s,t,u)\,dt\wedge du + q(s,t,u)\,du\wedge ds + r(s,t,u)\,ds\wedge dt\,.$$
Then we have
\begin{align*}
	p(s, t, u) & = F^\dast\overline{\omega}(\partial_t,\partial_u)
=\overline{\omega}(F_\dast\partial_t,F_\dast\partial_u)\\
	& =\overline{\omega}(F_t,F_u)=\frac{(F_t\times F_u)\dtw F}{4\pi|F|^3}\,,
\end{align*}
where the subscripts on $F$ denote partial derivatives, and likewise for 
$q(s, t, u)$ and $r(s, t, u)$.

Therefore, the characteristic 2-form of the link $L$ is 
\begin{align*}
	\omega_L &= g_L^\dast\omega \\
	&=\frac{1}{4\pi|F|^3}\left(\vphantom{\sqrt{a^2}}(F_t\times F_u\dtw F)\,dt\wedge du
+(F_u\times F_s\dtw F)\,du\wedge ds+(F_s\times F_t\dtw F)\,ds\wedge dt\right)\,.
\end{align*}

%%%%%%%%%%%%%%%%%%%%%%%%%%%%%
%%%%%%%%%%%%%%%%%%%%%%%%%%%%%
\part{Proof of \thmref{B}, formula \eqref{eq:B1}.}
%%%%%%%%%%%%%%%%%%%%%%%%%%%%%
%%%%%%%%%%%%%%%%%%%%%%%%%%%%%

Let $L$ be a three-component link in Euclidean 3-space $\reals^3$ with pairwise linking numbers $p$, $q$ and $r$ all zero.

By the first part of \thmref{A} these numbers are the degrees of the Gauss map $g_L\colon T^3\to S^2$ on the two-dimensional coordinate subtori. Since these degrees are all zero, $g_L$ is homotopic to a map $g\colon T^3\to S^2$ which collapses the 2-skeleton of $T^3$ to a point:
\begin{diagram}[size=2.2em]
	g_L  \simeq g\colon & T^3 & \rTo^{\sigma}& S^3 & \rTo^{f} & S^2
\end{diagram}
where $\sigma$ is the collapsing map. 

By the second part of \thmref{A}, Milnor's $\mu$-invariant of $L$ is equal to half of Pontryagin's $\nu$-invariant of $g_L$, which in turn is just the Hopf invariant of $f\colon S^3\to S^2$,
$$\mu(L) = \frac12 \nu(g_L) = \frac12\Hopf(f)\,.$$ 
We can thus use J.H.C.\ Whitehead's integral formula for the Hopf invariant, as follows.

Let $\omega$ be the area 2-form on $S^2$, normalized so that $\int_{S^2}\omega=1$. 
Its pullback $f^\dast\omega$ is a closed 2-form on $S^3$, which is exact 
because $H^2(S^3;\reals) = 0$. Let $d^{-1}(f^\dast\omega)$ indicate any smooth 
1-form on $S^3$ whose exterior derivative is $f^\dast\omega$. Then, as noted 
earlier, Whitehead showed that the Hopf invariant of f is given by the formula
$$\Hopf(f) = \int_{S^3}d^{-1}(f^\dast\omega)\wedge f^\dast\omega\,,$$
the value of the integral being independent of the choice of the 1-form 
$d^{-1}(f^\dast\omega)$.

Pulling the integral back to $T^3$ via the collapsing map $\sigma\colon T^3\to S^3$
yields the formula 
\[
	\nu(g_L) = \int_{T^3} d^{-1}(\omega_L)\wedge \omega_L\, ,
\]
thanks to the fact that $g_L$ is homotopic to $g = f\circ\sigma$, 
and recalling that
$\omega_L=g_L^\dast\omega$. Since $\mu(L) = \frac12 \nu(g_L)$, we get 
$$\mu(L) = \frac12 \int_{T^3} d^{-1}(\omega_L)\wedge\omega_L\,,$$
completing the proof of formula \eqref{eq:B1} of \thmref{B}.

%%%%%%%%%%%%%%%%%%%%%%%%%%%%%
%%%%%%%%%%%%%%%%%%%%%%%%%%%%%
\part{The geometrically natural choice for $d^{-1}(\omega_L)$.}
%%%%%%%%%%%%%%%%%%%%%%%%%%%%%
%%%%%%%%%%%%%%%%%%%%%%%%%%%%%

After stating \thmref{B}, we indicated that the geometrically natural choice for 
$d^{-1}(\omega_L)$ 
is the primitive of least $L^2$-norm. It can be obtained explicitly by convolving 
$\omega_L$ with the fundamental solution $\varphi$ of the scalar Laplacian on $T^3$, 
and then taking the exterior co-derivative $\delta$ of the resulting 2-form:
$$d^{-1}(\omega_L)=\delta(\varphi\ast\omega_L)\,.$$
If we make this geometrically natural choice for $d^{-1}(\omega_L)$, then the integrand
$$d^{-1}(\omega_L)\wedge\omega_L$$
in the above formula for the triple linking number $\mu(L)$ is also geometrically 
natural in the sense that it is unchanged if $L$ is moved by an orientation-preserving 
isometry of $\reals^3$.

We give a hint of the details, extracted from Part 1, and refer the reader there for 
proofs of Propositions~\ref{prop:A} and~\ref{prop:B} below.

\begin{proposition}\label{prop:A}
	\textbf{\boldmath The fundamental solution 
of the scalar Laplacian on the 3-torus ${T^3 =(\reals/ 2\pi\ints)^3}$ is given by the formula
	$$\varphi(\x) = \frac{1}{8\pi^3} \sum_{\n\neq\zero} \frac{e^{i\n\dte\x}}{|\n|^2}\,.$$ 
	The function $\varphi$ is $C^\infty$ at all points $\x\in T^3$ except $\zero$, where it becomes infinite.}
\end{proposition}

In the above formula, $\n$ denotes a triple of integers.

\begin{figure}[h!]
	\includegraphics[width=230pt]{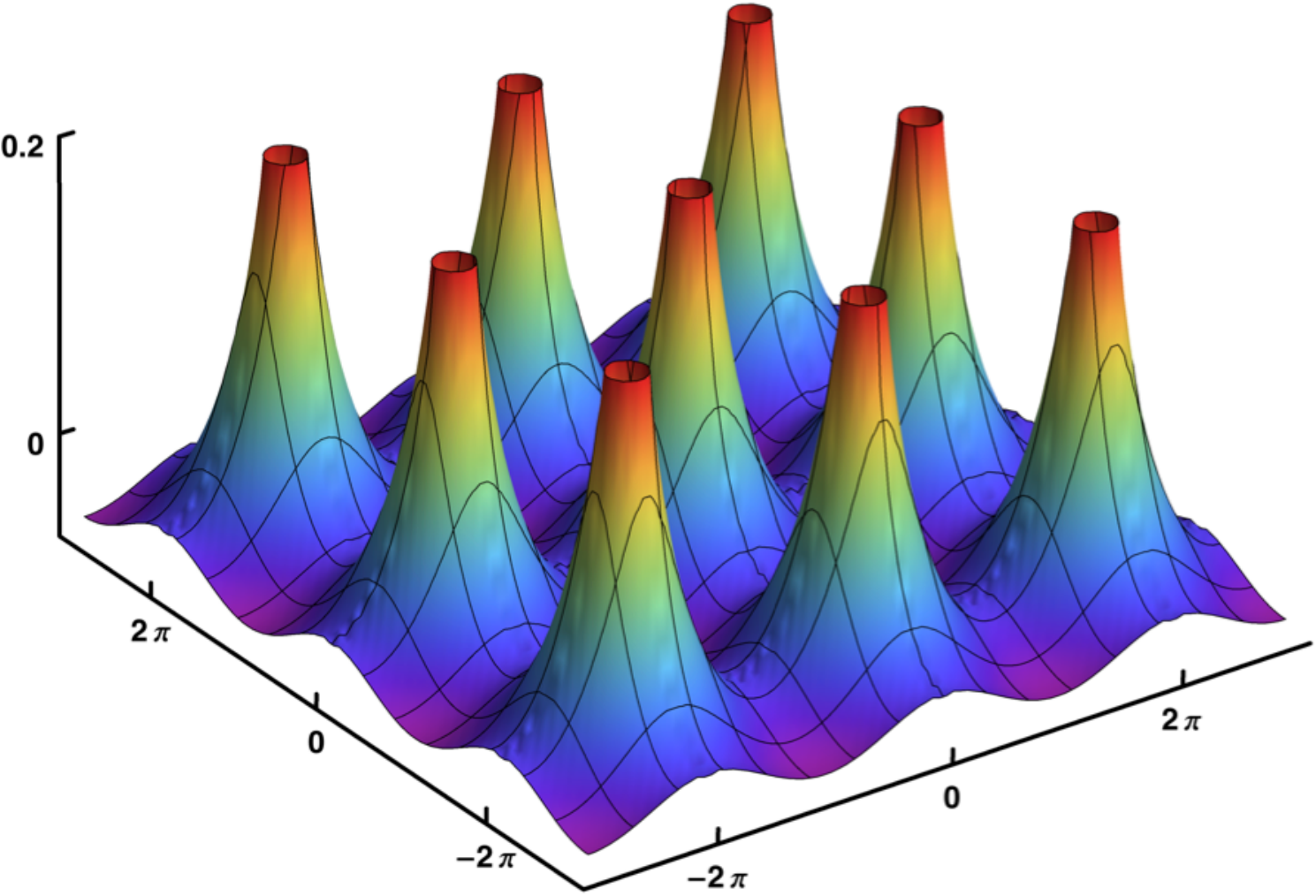}
	\caption{Fundamental solution of the Laplacian}
	\label{fig:fundsol}
\end{figure}

Although this formula for $\varphi$ is expressed in terms of complex exponentials, the value of $\varphi$ is real for real values of $\x$ because of the symmetry of the coefficients. Figure~\ref{fig:fundsol} shows the graph of the corresponding fundamental solution 
$$\varphi(\x) = \frac{1}{4\pi^2}\sum_{\n\neq\zero}\frac{e^{i\n\dt\x}}{|\n|^2}$$
of the scalar Laplacian on the 2-torus $T^2	= S^1\times S^1$, summed for $|\n|\le 15$, and displayed over the range $[-3\pi, 3\pi]\times [-3\pi, 3\pi]$.

%%%% FIGURE 12 %%%%%
%%%%%%%%%%%%%%%%%%%

The graph in Figure~\ref{fig:fundsol} looks like a ``Morse function'' with infinite maxes at the lattice points, saddles in the middle of the ``edges'', and mins at the center of the fundamental domains.

\begin{proposition}\label{prop:B}
	\textbf{\boldmath If $\omega$ is any exact 
differential form on $T^3$ with $C^\infty$ coefficients, then
$$\alpha=\delta(\varphi\ast\omega)$$
is a $C^\infty$ differential form satisfying $d\alpha=\omega$. Furthemore, 
if $d\overline{\alpha}=\omega$
as well, then $|\alpha|_{L^2}\le|\overline{\alpha}|_{L^2}$, with equality if and only 
if $\overline{\alpha}=\alpha$.}
\end{proposition}

%%%%%%%% THIS IS THE END OF SECTION 3 %%%%%%%%

%% file: 4epilogue.tex
%%%%%%%%%%%%%%%%%%%%%%%%%
%% EPILOGUE
%%%%%%%%%%%%%%%%%%%%%%%%%

\clearpage

\section*{Epilogue} 
\label{sec:epilogue}

%%%%%%%%%%%%%%%%%%%%%%%%%%%%%
%%%%%%%%%%%%%%%%%%%%%%%%%%%%%
\part{Where do the generalized Gauss maps come from?}
%%%%%%%%%%%%%%%%%%%%%%%%%%%%%
%%%%%%%%%%%%%%%%%%%%%%%%%%%%%

In the spherical theory, the generalized Gauss map $g_L\colon T^3\to S^2$ comes from the key map
\begin{diagram}
	f_S: & \confs & \rTo^{\text{Grassmann map } G} & G_2 \br^4 & \rTo^{\pi_+} & S^2
\end{diagram}	
via the substitution
$$g_L(s, t, u) = f_S\left(x(s) , y(t) , z(u)\right)\,,$$ 
while in the Euclidean theory it comes in the same way from the unit normalization $f_E$ of the key map 
$$F_E(x,y,z) = [a]+[b]+[c]+[b,c]+[c,a]+[a,b]\,.$$
But where do these key maps come from?

In the spherical theory, we saw in Part 1 that the configuration space $\confs$ deformation retracts to a subspace diffeomorphic to $S^3 \times S^2$, and the key map $f_S$ there is an $SO(4)$-equivariant version of this deformation retraction, followed by projection to the $S^2$ factor.

In the Euclidean theory, we face two complicating features: the configuration space $\conf$ is more challenging -- it deformation retracts to a subspace diffeomorphic to a nontrivial $S^2 \vee S^2$ bundle over $S^2$ -- and the group $\Isom^+\reals^3$ of orientation-preserving isometries of $\reals^3$ is non-compact.

If we were not seeking a generalized Gauss map which is geometrically natural in the sense of being $\Isom^+\reals^3$-equivariant, we could simply define the key map $f_E$ to be the composition
\begin{diagram}
	\conf & \rTo[l>=.5in]^H & \confs & \rTo[l>=.5in]^{f_S}&  S^2\,,
\end{diagram}
where the ``inclusion'' $H\colon \conf \longrightarrow \confs$ was defined earlier via inverse stereographic projection. This definition of $f_E$ is far from being $\Isom^+\reals^3$-equivariant, and the resulting generalized Gauss map would suffer from the same defect, and so lose its applicability to problems in fluid dynamics and plasma physics.

What we did instead was to consider the map $H'\colon \conf \longrightarrow \confs$ which first took a triple $(x, y, z)$ of distinct points in $\reals^3$ via translation to the triple $(x-z,y-z,0)$ ``based'' at the origin in $\reals^3$, and then via inverse stereographic projection $h$ to the triple $(h(x-z) , h(y-z) , 1)$ of distinct points based at the identity in $S^3$.

Inverse stereographic projection for such based triples is $SO(3)$-equivariant, and leads to the map $F_S\circ H'\colon\conf \longrightarrow \reals^3-\{0\}$ which up to scale takes
\begin{diagram}
	(x,y,z) & \rMapsto & |a|^2 b + a|b|^2 + a\times b\,,
\end{diagram}
thanks to our earlier computation, where $a = z-y$ and $b = x-z$. This map is $\Isom^+\reals^3$-equivariant, and so is the resulting projection to $S^2$, but at the cost of losing scale-invariance and ``sign symmetry'' in the three points $x$, $y$ and $z$.

A little artful play led to the alternative formula 
$$F_E(x,y,z) =[a]+[b]+[c]+[b,c]+[c,a]+[a,b]\,,$$
which is still $\Isom^+\reals^3$-equivariant, but now also scale-invariant and sign symmetric, and at the same time, thanks to the Bridge Lemma, homotopic to $F_S\circ H'$.

This is the origin of the key map $F_E\colon\conf\to\reals^3-\{0\}$ and the resulting Euclidean version of the generalized Gauss map $g_L\colon T^3\to S^2$.

%%%%%%%% THIS IS THE END OF THE EPILOGUE %%%%%%%%

%% file: 5bibliography.tex
%%%%%%%%%%%%%%%%%
%% REFERENCES
%%%%%%%%%%%%%%%%%

\providecommand{\bysame}{\leavevmode\hbox to3em{\hrulefill}\thinspace}
\providecommand{\MR}{\relax\ifhmode\unskip\space\fi MR }
\providecommand{\MRhref}[2]{%
  \href{http://www.ams.org/mathscinet-getitem?mr=#1}{#2}
}
\providecommand{\href}[2]{#2}

\bigskip

\noindent
\it deturck@math.upenn.edu\rm\\
\it gluck@math.upenn.edu\rm\\
\it rako@tulane.edu\rm\\ 
\it pmelvin@brynmawr.edu\rm\\ 
\it hnuchi@math.upenn.edu\rm\\ 
\it clayton@math.uga.edu\rm\\ 
\it shea@math.lsu.edu\rm\\